\documentclass[12pt,draft]{article}
\usepackage{amsmath,amsfonts,amssymb,amsthm,amscd,comment}

\newcommand{\w}{\omega}
\newcommand{\1}{{\bf 1}}
\newcommand{\id}{{\rm id}}

\DeclareMathOperator{\Ker}{Ker\,}

\newcommand{\wt}[1]{|#1|}

\newcommand{\image}{{\rm Im\,}}
\newcommand\hari[2]{\left\langle\,#1\,\left\vert \,#2\,\right.\right\rangle_{\C}}
\newcommand\haru[2]{\left\langle\left.\,#1\,\right\vert \,#2\,\right\rangle_{\C}}

\newcommand\Z{\mathbb{Z}}
\newcommand\Zpos{\Z_{\geq0}}
\newcommand\Zplus{\Z_{>0}}
\newcommand\Q{\mathbb{Q}}

\newcommand\C{\mathbb{C}}

\newcommand{\NO}{\,{\raise0.25em\hbox{$\mathop{\hphantom {\cdot}}\limits^{_{\circ}}_{^{\circ}}$}}\,}

\newtheorem{theorem}{Theorem}[section]
\newtheorem{proposition}[theorem]{Proposition}
\newtheorem{lemma}[theorem]{Lemma}
\newtheorem{corollary}[theorem]{Corollary}

\theoremstyle{definition}

\theoremstyle{remark}

\numberwithin{equation}{section}

\begin{document}
\begin{center}
\begin{Large}
$C_2$-cofiniteness of the $2$-cycle permutation orbifold models of minimal Virasoro vertex operator algebras
\end{Large}
\vskip1ex
Toshiyuki Abe\footnote{This work was supported by Grant-in-Aid for Young Scientists (B) 20740017}\\
Graduate school of Science and Engineering, Ehime university\\
2-5, Bunkyocho, Matsuyama, Ehime 790-8577, Japan
 
\vskip1cm
{\bf Abstract}
\end{center}
\begin{small}
In this article, we give a sufficient and necessary condition for the $C_2$-cofiniteness of $\widetilde{V}=(V\otimes V)^\sigma$ for a $C_2$-cofinite vertex operator algebra $V$ and the $2$-cycle permutation $\sigma$ of $V\otimes V$.  
As an application, we show that the $2$-cycle permutation orbifold model of the simple Virasoro vertex operator algebra $L(c,0)$ of minimal central charge $c$ is $C_2$-cofinite.
\end{small}

\section{Introduction}
Given a vertex operator algebra $V$, the tensor product $V^{\otimes k}$ of $k$-copies of $V$ as a $\C$-vector space canonically has a vertex operator algebra structure (see \cite{FHL}).
Each permutation $\sigma$ of tensor factors gives rise to an automorphism of $V^{\otimes k}$ of finite order, and the fixed point set $(V^{\otimes k})^\sigma=\{u\in V^{\otimes k}\,|\,\sigma(u)=u\}$ becomes a vertex operator subalgebra called a $\sigma$-permutation orbifold model.
In this article, we consider the $2$-transposition $\sigma$ of $V\otimes V$ for a $C_2$-cofinite vertex operator algebra $V$, and study the $C_2$-cofiniteness condition of $(V\otimes V)^\sigma$.
We give a sufficient and necessary condition for the $C_2$-cofiniteness of the $\sigma$-permutation orbifold model, and we show that $(L(c,0)\otimes L(c,0))^\sigma$ is $C_2$-cofinite for any minimal central charge $c=c_{p,q}=1-6\frac{(p-q)^2}{pq}$ with coprime integers $p,q\geq 2$.

The $C_2$-cofiniteness condition requires that $V/\haru{a_{(-2)}b}{a,b\in V }$ is finite dimensional, where $a_{(-2)}b$ denotes the $-2$-product of ordered pair $(a,b)$ for $a,b\in V$.
A vertex operator algebra satisfying this condition is called {\it $C_2$-cofinite}.
The $C_2$-cofiniteness condition is one of the most important properties in the vertex operator algebra theory, and once a vertex operator algebra satisfies this condition, its modules have a lot of remarkable features, for example, modular invariance of (extended) trace functions and the existence of fusion products (see \cite{Huang07}, \cite{Miyamoto04}--\cite{Miyamoto10} and references therein).
In spite of its importance, a verification of the $C_2$-cofiniteness condition is a very difficult task in general, and some conjectures for the $C_2$-cofiniteness have been still remained open.
One of them is that if $V$ is a $C_2$-cofinite simple vertex operator algebra and $G$ is a finite automorphism group, then the fixed point vertex operator subalgebra $V^G:=\{a\in V\,|\,g(a)=a,\forall g\in G\}$ is $C_2$-cofinite.
This conjecture has been checked only for some examples such as lattice vertex operator algebras with a certain involution (\cite{Yamskulna09}), and no definitive idea to prove this conjecture has been found.  

The permutation orbifold theory has been studied by physicists two decades ago(see \cite{KlemmSchmidt90}, \cite{FuchsKlemmSchmidt92}, see also \cite{Bantay02}).
Its systematic study has been started with the papers \cite{BorisovHalpernSchweigert98} for cyclic permutations, and was generalized to the general symmetric group in \cite{Bantay98}.
The techniques are effectively applied to the study of elliptic genera in terms of the second quantized string theory in \cite{DijkgraafMooreVerlindeVerlinde97} and also applied to the study of the kernels of the modular representations in \cite{Bantay03} for example.   
On the other hand, the remarkable paper \cite{BarronDongMason02} seems to be known as a unique paper in which permutation orbifold models are studied from the point of view of the vertex operator algebra theory.
The aim of this article is to start the study of structure theory of permutation orbifold models in the theory of vertex operator algebras, and especially to show the $C_2$-cofiniteness condition is valid for permutation orbifold models of $C_2$-cofinite vertex operator algebras.
As the very first step, we consider $2$-cycle permutation orbifold models of $C_2$-cofinite vertex operator algebras.
We could not show that arbitrary permutation orbifold model is $C_2$-cofinite without any condition except a simpleness and $C_2$-cofiniteness of the based vertex operator algebra.
But we give a sufficient and necessary condition for this assertion.

Let $V$ be a vertex operator algebra.  
Then the $2$-cycle permutation orbifold model $\widetilde{V}=(V\otimes V)^\sigma$ is linearly spanned by vectors $\phi(a,b)=a\otimes b+b\otimes a$ for $a,b\in V$.
We then have an identity
\begin{align}\label{gdfuuq}
\phi(a_{(-n)}u,v)\equiv -\phi(u,a_{(-n)}v)\mod C_2(\widetilde{V})
\end{align}
for $a,u,v\in V$ and $n\geq 2$, where $a_{(-n)}b$ denotes the $-n$-product.
This identity plays an important role in our arguments.
Suppose that $V$ is $C_2$-cofinite and assume that $V$ is strongly generated by a finite set $T$.
Then by using the identity \eqref{gdfuuq} and other facts, we can prove that if $V$ is $C_2$-cofinite, then the finiteness of the dimension of the subspace
\[
\haru{\phi(a_{(-n)}b,\1)+C_2(\widetilde{V})}{a,b\in T,n>0}
\]
 in $\widetilde{V}/C_2(\widetilde{V})$ is equivalent to the $C_2$-cofinitenss of $\widetilde{V}$.

For the case $V$ is the simple Virasoro vertex operator algebra $L(c,0)$ of a central charge $c$, $V$ is strongly generated by the Virasoro vector $\w$.
It is also known that $V$ is $C_2$-cofinite if and only if $c$ is a minimal central charge $c_{p,q}$ with relatively prime integers $p,q\geq 2$ (see \cite{Arakawa10} for example).
Therefore under the assumption that $c$ is a minimal central charge, if the subspace of $\widetilde{V}/C_2(\widetilde{V})$ spanned by $\phi(\w_{(-n)}\w,\1)+C_2(\widetilde{V})$ for $n\in\Zplus$ is finite dimensional then $\widetilde{V}$ is $C_2$-cofinite.
We in fact show that
\begin{align}\label{auuesss}
\phi(\w_{(-n)}\w,\1) \in C_2(\widetilde{V})
\end{align} for $n\geq 30$ regardless of the central charge $c$ by direct computations.
As a result we find that $\widetilde{V}$ is $C_2$-cofinite in the case $c$ is a minimal central charge.

This article is organized as follows.
In Section \ref{sect2} we recall some notions and results from representation theory of vertex operator algebras, and give some identities which are used in the later sections.   
In Section \ref{sect3}, we prove  \eqref{gdfuuq} and some useful lemmas.
In Section \ref{sect4}, we give a sufficient and necessary condition for the $C_2$-cofiniteness of $\widetilde{V}$ when $V$ is $C_2$-cofinite.
Section \ref{sect5} is devoted to calculations concerning with the Virasoro vectors.
We have some numerical results with the help of a computer and show that $(L(c,0)\otimes L(c,0))^\sigma$ is $C_2$-cofinite for a minimal central charge $c$.
Section \ref{sect6} is an appendix where we give a list of results of computations which is necessary in Section \ref{sect5} and show that \eqref{auuesss} is valid for $n\geq 30$.   

Acknowledgment:
The author would like to thank Professor Masahiko Miyamoto for insightful comments and helpful discussion.
He is also thankful to Prof. Hiromich Yamada for his encouragement and finding some typo in the first version. 
He greatly appreciate the referees for giving him useful comments and suggestions.

\section{Preliminaries}\label{sect2}
A vertex operator algebra is a quadruplet $(V,\cdot_{(n)}\cdot,\1,\w)$ consisting of a $\Zpos$-graded $\C$-vector space $V=\bigoplus_{d=0}^\infty V_d$, bilinear maps $V\times V\ni (a,b)\mapsto a_{(n)}b\in V$ associated to each integer $n\in\Z$, and two distinguished vectors $\1\in V_0$ called the vacuum vector and $\w\in V_2$ called the Virasoro vector.
For its axioms, refer \cite{MatsuoNagatomo99} for example.
We write down some identities from the axiom of the vertex operator algebras which we need in this article:   

{\bf Associativity formula:}
\begin{align}\label{asso1}
(a_{(m)}b)_{(n)}u=\sum_{i=0}^\infty\binom{m}{i}(-1)^i(a_{(m-i)}b_{(n+i)}u-(-1)^mb_{(m+n-i)}a_{(i)}u)
\end{align}
for $a,b,u\in V$ and $m,n\in\Z$.

{\bf Commutativity formula:}
\begin{align}\label{comm1}
[a_{(m)},b_{(n)}]u=\sum_{i=0}^\infty\binom{m}{i}(a_{(i)}b)_{(m+n-i)}u
\end{align}
for $a,b,u\in V$ and $m,n\in\Z$.

The vacuum vector satisfies that $a_{(i)}\1=0$ for $i\in\Zpos$, $\1_{(n)}a=\delta_{n,-1}a$.
As for the Virasoro vector $\w$, if we set $L_n=\w_{(n+1)}$ then
\begin{align}
[L_{m},L_{n}]&=(m-n)L_{m+n}+\delta_{m+n,0}\binom{m+1}{3}\frac{c_V}{2}\id_{V},\\
(L_{-1}a)_{(n)}&=-na_{(n-1)},\\
L_0a&=da\quad\text{for }d\in\Zpos\text{ and }a\in V_d.
\end{align}
Through this article, we assume that $V$ is of CFT type, that is, $V_0=\C\1$.

We give an identity deduced from the associativity formula.
\begin{lemma}\label{aiuewy}
For $a,b,u\in V$ and $m,n\in\Zplus$,
\begin{align*}
(a_{(-m)}b_{(-n)}\1)_{(-1)}u&=a_{(-m)}b_{(-n)}u\\
&\quad+\sum_{i=0}^\infty (\alpha_{m,n;i} a_{(-m-n-i)}b_{(i)}u+\alpha_{n,m;i}b_{(-m-n-i)}a_{(i)}u),
\end{align*}
where
\begin{align}\label{jashasu}
\alpha_{m,n;i}=\frac{(m+n-1)!}{(m-1)!(n-1)!}\binom{m+n-1+i}{i}\frac{(-1)^{n-1}}{n+i}
\end{align}
for $m,n\in\Zplus$ and $i\in\Zpos$.
\end{lemma}
\begin{proof}
By Associativity formula, we have
\begin{align}
\begin{split}\label{phyrd}
&(a_{(-m)}b_{(-n)}\1)_{(-1)}u\\
&=\sum_{i=0}^\infty \binom{-m}{i}(-1)^ia_{(-m-i)}(b_{(-n)}\1)_{(-1-i)}u\\
&\quad-\sum_{i=0}^\infty \binom{-m}{i}(-1)^{i+m}(b_{(-n)}\1)_{(-m-1-i)}a_{(i)}u\\
&=\sum_{i=0}^\infty \binom{m-1+i}{m-1}\binom{-1+i}{n-1}(-1)^{n-1}a_{(-m-i)}b_{(-n+i)}u\\&\quad-\sum_{i=0}^\infty\binom{m-1+i}{m-1}\binom{-m-1-i}{n-1}(-1)^{m+n-1}b_{(-m-n-i)}a_{(i)}u.
\end{split}
\end{align}
Since $\binom{-1+i}{n-1}=0$ if $1\leq i\leq n-1$, the first summation in the right hand side is equal to
\begin{align*}
&a_{(-m)}b_{(-n)}u+\sum_{i=0}^\infty \binom{m+n-1+i}{m-1}\binom{n-1+i}{n-1}(-1)^{n-1}a_{(-m-n-i)}b_{(i)}u
\end{align*}
We see that  
\begin{align*}
\binom{m+n-1+i}{m-1}\binom{n-1+i}{n-1}(-1)^{n-1}=\frac{(m+n-1+i)!(-1)^{n-1}}{(m-1)!(n-1)!i!(n+i)}=\alpha_{m,n;i}.
\end{align*}
On the other hand, it follows from $\binom{-m-1-i}{n-1}=\binom{m+n-1+i}{n-1}(-1)^{n-1}$ that the second summation of the right hand side in \eqref{phyrd} is equal to
\begin{align*}
&\sum_{i=0}^\infty\binom{m-1+i}{m-1}\binom{m+n-1+i}{n-1}(-1)^{m-1}b_{(-m-n-1-i)}a_{(i)}u\\
&=\sum_{i=0}^\infty\alpha_{n,m;i}b_{(-m-n-1-i)}a_{(i)}u.
\end{align*}
Hence we have the lemma.
\end{proof}

We next recall the definition of the $C_2$-cofiniteness condition and related results following \cite{Zhu96}.
A vertex operator algebra $V$ is said to satisfy the $C_2$-cofiniteness condition or to be $C_2$-cofinite if the subspace
\[
C_2(V):=\haru{a_{(-2)}b}{a,b\in V}
\]
is finite codimensional in $V$.
We denote by $\overline{\cdot}:V\rightarrow V/C_2(V)$ the canonical projection.
The quotient space $V/C_2(V)$ becomes a Poisson algebra with multiplication and Lie bracket
\[
\overline{a}\cdot\overline{b}=\overline{a_{(-1)}b},\quad [\overline{a},\overline{b}]=\overline{a_{(0)}b}
\]  
for $a,b\in V$ respectively.
In this article we do not use its Lie algebra structure.

A vertex operator algebra $V$ is called {\it strongly generated} by a subset $T\subset V$ if $V$ is spanned by $\1$ and vectors of the form $a^1_{(-n_1)}\cdots a^r_{(-n_r)}\1$ with $r\geq 1$, $a^i\in T$ and $n_i\geq 1$.
\begin{proposition}\label{uuuwuwu}
If $V$ is strongly generated by a set $T$, then $V/C_2(V)$ is generated by $\overline{T}$ as a commutative associative algebra with unit $\overline{\1}$.
\end{proposition}
Proposition \ref{uuuwuwu} can be proved by the facts that for $a^i\in T$ and $n_i\in\Zplus$,
\[
\overline{a^1_{(-1)}\cdots a^r_{(-1)}\1}=\overline{a^1}\cdot\cdots\cdot\overline{a^r}.
\]
and 
\[
\overline{a^1_{(-n_1)}\cdots a^r_{(-n_r)}\1}=0,
\]
if some $n_i>1$.
The following theorem is one of the important properties of the subspace $C_2(V)$.
\begin{theorem}\label{uasishia} {\rm (\cite{GaberdielNeitzke03})}
Let $V$ be a vertex operator algebra of CFT type, and $S$ a subset of $V$ such that $V=\C\1\oplus\langle S\rangle_{\C}\oplus C_2(V)$.
Then $V$ is spanned by $\1$ and vectors of the form $a^1_{(-n_1)}\cdots a^r_{(-n_r)}\1$ satisfying $r\geq 1$, $a^i\in S$, $n_1>n_2>\cdots >n_r\geq 1$.
\end{theorem}
By Theorem \ref{uasishia}, any $C_2$-cofinite vertex operator algebra has a finite subset $T$ by which $V$ is strongly generated.

\section{$2$-cyclic permutation orbifold models}\label{sect3}
Let $V$ be a vertex operator algebra.
Then $V\otimes V$ becomes a vertex operator algebra with Virasoro vector $\w\otimes\1+\1\otimes\w$ and vacuum vector $\1\otimes \1$.
The $n$-th product of $a\otimes b$ and $u\otimes v$ for $a,b,u,v\in V$, $n\in\Z$ is given by
\[
(a\otimes b)_{(n)}(u\otimes v)=\sum_{i\in\Z}(a_{(i)}u)\otimes (b_{(n-1-i)}v).
\]
In particular we have $C_2(V\otimes V)=C_2(V)\otimes V+V\otimes C_2(V)$.
Thus if $V$ is $C_2$-cofinite then so is $V\otimes V$.

We consider a map $\sigma\in GL(V\otimes V)$ defined by $\sigma(a\otimes b)=b\otimes a$ for $a,b\in V$.
Then we see that $\sigma$ is an automorphism of $V\otimes V$ of order $2$.  
Now we set
\[
\widetilde{V}:=(V\otimes V)^\sigma.
\]
Then $\widetilde{V}$ is a vertex operator subalgebra of $V\otimes V$ with vacuum vector $\1\otimes \1$ and Virasoro vector $\w\otimes\1+\1\otimes\w$.

We introduce some notations.  
For $a,b\in V$, we set   
\begin{align*}
\phi(a, b)=a\otimes b+b\otimes a,\quad\text{and }\eta(a)&=a\otimes \1+\1\otimes a\in\widetilde{V}.
\end{align*}
It is clear that $\phi(a,b)=\phi(b,a)$ for $a,b\in V$ and $\eta(a)=\phi(a,\1)=\phi(\1,a)$.
\begin{lemma}\label{aisu}
For any $a,u,v\in V$ and $n\in\Z$,
\begin{align}\label{uygfqw}
\eta(a)_{(n)}\phi(u,v)=\phi(a_{(n)}u,v)+\phi(u,a_{(n)}v).
\end{align}
\end{lemma}
\begin{proof}
We have
\begin{align*}
\eta(a)_{(n)}\phi(u,v)&=(a\otimes \1+\1\otimes a)_{(n)}(u\otimes v+v\otimes u)\\
&=a_{(n)}u\otimes v+a_{(n)}v\otimes u+u\otimes a_{(n)}v+b\otimes a_{(n)}v\\
&=\phi(a_{(n)}u,v)+\phi(u,a_{(n)}v).
\end{align*}
This proves the lemma.  
\end{proof}

\begin{lemma}\label{uhais}
For any $a,b\in V$, $\eta(a)_{(-1)}\eta(b)=\eta(a_{(-1)}b)+\phi(a,b)$.
\end{lemma}
\begin{proof}
Substitute $u=b$, $v=\1$ and $n=-1$ in \eqref{uygfqw}.  
\end{proof}
Since $\widetilde{V}$ is spanned by $\image\phi$, we have the following proposition.
\begin{proposition}
$\widetilde{V}$ is strongly generated by $\image \eta$.
\end{proposition}
Now we consider $\widetilde{V}/C_2(\widetilde{V})$.
We denote by $\overline{\phi}$ and $\overline{\eta}$ the compositions of $\eta$ and $\phi$ with the canonical projection $\overline{\cdot}:\widetilde{V}\rightarrow \widetilde{V}/C_2(\widetilde{V})$ respectively.  
Since $L_{-1}(\widetilde{V})\subset C_2(\widetilde{V})$, we have
\begin{align}\label{dyeh}
\overline{\phi}(L_{-1}a,b)=-\overline{\phi}(a,L_{-1}b)
\end{align}
for any $a,b\in V$.
In particular one sees that $\overline{\eta}(L_{-1}a)=0$.
Therefore we have
\begin{align}\label{iwweyyo}
L_{-1}V\subset \Ker \overline{\eta},
\end{align}
and
\begin{align}\label{aiaio}
\overline{\eta}(a_{(-n)}\1)=0
\end{align}
for $a\in V$ and $n\geq 2$.

The following lemma is an important property of $\overline{\phi}$.
\begin{lemma}\label{apjrj}
For any $a,b,u\in V$ and $n\geq 2$, $\overline{\phi}(a_{(-n)}u,b)=-\overline{\phi}(u,a_{(-n)}b) $.
\end{lemma}
\begin{proof}
By Lemma \ref{aisu},
\begin{align*}
\phi(a_{(-n)}u,b)+\phi(u,a_{(-n)}b)=\eta(a)_{(-n)}\phi(u,b)\in C_2(\widetilde{V}).
\end{align*}
Thus we have the lemma.
\end{proof}
Now we show more complicated identities.
\begin{lemma}\label{asyudgf}
For $a,b,u\in V$ and $n\geq 2$,
\begin{align}
\begin{split}\label{uyagsu}
&\overline{\phi}(a,b_{(-n)}u)+\overline{\phi}(b,a_{(-n)}u)\\
&=\sum_{i=0}^{\infty}\left(\overline{\eta}(a_{(-i-2)}b_{(-n+1+i)}u)+\overline{\eta}(b_{(-i-2)}a_{(-n+1+i)}u)\right).
\end{split}
\end{align}
\end{lemma}
\begin{proof}
We have
\begin{align*}
\phi(a,b)_{(-n)}\eta(u)&=(a\otimes b+b\otimes a)_{(-n)}(u\otimes \1+\1\otimes u)\\
&=\sum_{i\in\Z}(a_{(i)}u\otimes b_{(-n-1-i)}\1+a_{(i)}\1\otimes b_{(-n-1-i)}u\\
&\quad +b_{(i)}u\otimes a_{(-n-1-i)}\1+b_{(i)}\1\otimes a_{(-n-1-i)}u)\\
&=\sum_{i=-n}^\infty (a_{(i)}u\otimes b_{(-n-1-i)}\1+b_{(i)}u\otimes a_{(-n-1-i)}\1)\\
&\quad +\sum_{i=0}^\infty(a_{(-i-1)}\1\otimes b_{(-n+i)}u +b_{(-i-1)}\1\otimes a_{(-n+i)}u)\\
&=\sum_{i=0}^\infty (\phi(a_{(-n+i)}u,b_{(-1-i)}\1)+\phi(b_{(-n+i)}u, a_{(-1-i)}\1)).
\end{align*}
By Lemma \ref{apjrj}, we see that if $i>0$, then
\begin{align*}
&\overline{\phi}(a_{(-n+i)}u,b_{(-1-i)}\1)+\overline{\phi}(b_{(-n+i)}u, a_{(-1-i)}\1)\\
&= -\overline{\phi}(b_{(-1-i)}a_{(-n+i)}u,\1)-\overline{\phi}(a_{(-1-i)}b_{(-n+i)}u, \1)\\
&=-\sum_{i=0}^{\infty}\left(\overline{\eta}(a_{(-i-2)}b_{(-n+1+i)}u)+\overline{\eta}(b_{(-i-2)}a_{(-n+1+i)}u)\right).
\end{align*}
Since $\phi(a,b)_{(-n)}\eta(u)\in C_2(\widetilde{V})$, we have \eqref{uyagsu}.
\end{proof}

Therefore we have  
\begin{align}\label{asiw}
\overline{\phi}(a,b_{(-n)}u)= -\overline{\phi}(b,a_{(-n)}u)\mod \image \overline{\eta}
\end{align}
for $a,b,u\in V$ and $n\geq 2$.
Finally we show the following.
\begin{lemma}\label{uehwi}
Let $a,b,u\in V$ and $n\in\Zplus$.
If $n\geq 3$ then $\overline{\phi}(a,b_{(-n)}u)\in\image\overline{\eta}$.
\end{lemma}
\begin{proof}
Since $n-1\geq 2$, by using \eqref{uyagsu}, one has
\begin{align*}
\overline{\phi}(a,b_{(-n)}u)&=\frac{1}{n-1}\overline{\phi}(a, (L_{-1}b)_{(-n+1)}u)\\
&= -\frac{1}{n-1}\overline{\phi}(L_{-1}b, a_{(-n+1)}u)\mod \image \overline{\eta}\\
&= \frac{1}{n-1}\overline{\eta}(b_{(-2)}a_{(-n+1)}u)\in \image \overline{\eta}.
\end{align*}
\end{proof}

Suppose that $V$ is of CFT type and take a subset $S$ of $V$ such that
\[
V=\C\1\oplus\langle S\rangle_{\C}\oplus C_2(V).
\]
Then by Theorem \ref{uasishia}, $V$ is spanned by vectors of the form
\begin{align}\label{eydyu}
v={x^{1}}_{(-n_{1})}\cdots {x^{r}}_{(-n_{r})}\1,\quad x^j\in S,\quad n_{1}>\cdots >n_{r}\geq 1.
\end{align}
By \eqref{aiaio}, if $r=1$ and $n_1\geq 2$, then $\overline{\eta}(v)=0$.
Next we consider the case $r\geq 2$.
Since $n_1>\cdots >n_{r-1}\geq 2$, for $u\in V$, we have  
\[
\overline{\phi}(v,u)=(-1)^{r-1}\overline{\phi}({x^{r}}_{(-n_{r})}\1,{x^{r-1}}_{(-n_{r-1})}\cdots {x^{1}}_{(-n_{1})}u).
\]
If $n_r\geq 2$, then $\overline{\phi}(v,u)\in\image \overline{\eta}$ by Lemma \ref{apjrj}.
On the other hand in the case $n_r=1$, we can write ${x^{r-1}}_{(-n_{r-1})}\cdots {x^{1}}_{(-n_{1})}u$ as a linear combination of vectors of the form \eqref{eydyu}, say
\begin{align*}
w={y^{1}}_{(-m_{1})}\cdots {y^{r}}_{(-m_{s})}\1,\quad y^j\in S,\quad m_{1}>\cdots >m_{s}\geq 1.
\end{align*}
Then if $m_1\geq 3$ then $\overline{\phi}(x^{r},w)\in \image\overline{\eta}$ by Lemma \ref{uehwi}.
Consequently, we find that for any vector $v$ of the form \eqref{eydyu} and $u\in V$, $ \overline{\phi}(v,u)$ is equivalent to a linear combination of vectors of the form $\overline{\phi}(x,y)$, $\overline{\phi}(x,y_{(-1)}z)$ or $\overline{\phi}(x,y_{(-2)}z)$ for $x,y,z\in S$ modulo $\image \overline{\eta}$.  

Therefore we have
\begin{align}
\begin{split}\label{isauhq}
\widetilde{V}/C_{2}(\widetilde{V})=&\haru{\overline{\phi}(x,y)}{x,y\in S}\\
&+\haru{\overline{\phi}(x,y_{(-i)}z)}{x,y,z\in S, i=1,2}+\image \overline{\eta}.
\end{split}
\end{align}
Consequently we have the following lemma.
\begin{lemma}
If $V$ is $C_2$-cofinite and of CFT type, then the quotient space of $\widetilde{V}/C_2(\widetilde{V})$ by $\image\overline{\eta}$ is finite dimensional.
\end{lemma}
\begin{proof}
If $V$ is $C_2$-cofinite, then we may take $S$ to be a finite set.
Thus \eqref{isauhq} proves the lemma.
\end{proof}
Now the following corollary holds.
\begin{corollary}\label{isuhds}
Let $V$ be a $C_2$-cofinite vertex operator algebra of CFT type.
If $\image\overline{\eta}$ is finite dimensional, then $\widetilde{V}$ is $C_2$-cofinite.
\end{corollary}

\section{On the finiteness of $\dim \image \overline{\eta}$}\label{sect4}
In this section we consider the finiteness of dimension of the image $\image \overline{\eta}$ and show the following theorem:
\begin{theorem}\label{fbeygfw}
Let $V$ be a $C_2$-cofinite simple vertex operator algebra of CFT type, and suppose that $V$ is strongly generated by a finite subset $T\subset \cup_{d=1}^\infty V_d$.
Then $\widetilde{V}$ is $C_2$-cofinite if and only if $\haru{\overline{\eta}(x_{(-n)}y)}{x,y\in T, n\geq 1}$ is finite dimensional.
\end{theorem}
In Theorem \ref{fbeygfw}, the ``only if'' part is clear.
Hence it suffices to prove that the following lemma by Corollary \ref{isuhds}.
\begin{lemma}\label{ashgf}
Let $V$ and $T$ be as in Theorem \ref{fbeygfw}.
Then $\image\overline{\eta}$ is finite dimensional if $\haru{\overline{\eta}(x_{(-n)}y)}{x,y\in T, n\geq 1}$ is finite dimensional.  
\end{lemma}
We give a proof of Lemma \ref{ashgf} at the end of this section.  

We first recall the standard filtration of a vertex operator algebra (see \cite{GaberdielNeitzke03} or \cite{Miyamoto04}).
Set
\[
\mathcal{F}_{k}V:=\hari{a^1_{(n_1)}\cdots a^r_{(n_r)}\1}{r\geq 1, a^i\in V,\text{ homogeneous}, \sum_{j=1}^{r}\wt{a^j}\leq k}
\]
for $k\in\Zpos$.
The following proposition is well known (see \cite{Li04} for its proof).
\begin{proposition}\label{euysy}
Let $a^1,\cdots,a^r\in V$ be homogeneous vectors, $n_1,\cdots,n_r\in\Z$, $k,l\in\Zpos$ and $u\in\mathcal{F}_l(V)$.
\begin{enumerate}
\item[{\rm (1)}] If $\sum \wt{a^i}=k$, then ${a^1}_{(n_1)}\cdots {a^r}_{(n_r)}u\equiv {a^{i_1}}_{(n_{i_1})}\cdots {a^{i_r}}_{(n_{i_r})}u$ modulo the subspace $\mathcal{F}_{k+l-1}(V)$ for any permutation $\{i_1,\cdots,i_r\}$ of $\{1,\cdots,r\}$.  
\item[{\rm (2)}] If $\sum \wt{a^i}=k$ and $n_i\geq 0$ for some $i$, then  ${a^1}_{(n_1)}\cdots {a^r}_{(n_r)}u\in \mathcal{F}_{k+l-1}(V)$.
\item[{\rm (3)}] If $\sum \wt{a^i}=k$, then $({a^1}_{(n_1)}\cdots {a^r}_{(n_r)}\1)_{(-1)}u\equiv {a^1}_{(n_1)}\cdots {a^r}_{(n_r)}u$ modulo the subspace $\mathcal{F}_{k+l-1}(V)$.
\end{enumerate}
\end{proposition}

Now we consider a filtration $\overline{\eta}(\mathcal{F}_{k}V)$ for $k\in\Zpos$ of $\image \overline{\eta}\subset\widetilde{V}/C_2(\widetilde{V})$.
Then Proposition \ref{euysy} shows the following proposition.
\begin{proposition}\label{cdyuga}
Let $a^1,\cdots,a^r\in V$ be homogeneous, $n_1,\cdots,n_r\in\Z$, $k,l\in\Zpos$ and $u\in\mathcal{F}_l(V)$.
\begin{enumerate}
\item[{\rm (1)}] If $\sum \wt{a^i}=k$, then $\overline{\eta}({a^1}_{(n_1)}\cdots {a^r}_{(n_r)}u)\equiv \overline{\eta}({a^{i_1}}_{(n_{i_1})}\cdots {a^{i_r}}_{(n_{i_r})}u)$ modulo  $\overline{\eta}(\mathcal{F}_{k+l-1}(V))$ for any permutation $\{i_1,\cdots,i_r\}$ of $\{1,\cdots,r\}$.  
\item[{\rm (2)}] If $\sum \wt{a^i}=k$ and $n_i\geq 0$ for some $i$, then  $\overline{\eta}({a^1}_{(n_1)}\cdots {a^r}_{(n_r)}u)\equiv 0$ modulo $\overline{\eta}(\mathcal{F}_{k+l-1}(V))$.
\item[{\rm (3)}] If $\sum \wt{a^i}=k$, then $\overline{\eta}(({a^1}_{(n_1)}\cdots {a^r}_{(n_r)}\1)_{(-1)}u)\equiv \overline{\eta}({a^1}_{(n_1)}\cdots {a^r}_{(n_r)}u)$ modulo $\overline{\eta}(\mathcal{F}_{k+l-1}(V))$.
\end{enumerate}
\end{proposition}

Suppose that $V$ is strongly generated by a subset $T$ consisting of homogeneous vectors.
By definition, we see that $V$ is spanned by vectors of the form
\begin{align}\label{iyasd}
x^1_{(-n_1)}\cdots x^r_{(-n_r)}\1\quad\text{with }n_i\geq 1\text{ and }x_i\in T.
\end{align}
Moreover, each subspace $\mathcal{F}_{k}(V)$ is spanned by vectors of the form \eqref{iyasd} subject to $\sum_{i=1}^r\wt{x^i}\leq k$.

The following lemma plays an essential role in the proof of Lemma \ref{ashgf}.
\begin{lemma}\label{uaygsdf}
For any $a,b,u\in V$ and $m,n\geq 2$,  
\begin{align*}
&\overline{\eta}(a_{(-m)}b_{(-n)}\1)\cdot \overline{\eta}(u)\\
&=2\overline{\eta}(a_{(-m)}b_{(-n)}u)+\sum_{i=0}^\infty\binom{-n}{i}\overline{\eta}((b_{(i)}a)_{(-m-n-i)}u)\\
&\quad+\sum_{i=0}^\infty \left(\alpha_{m,n:i}\overline{\eta}(a_{(-m-n-i)}b_{(i)}u)+\alpha_{m,n;i}\overline{\eta}(b_{(-m-n-i)}a_{(i)}u)\right),
\end{align*}
where $\alpha_{m,n:i}$ is a constant defined in \eqref{jashasu}.
\end{lemma}
\begin{proof}
Set
\begin{align*}
w=\sum_{i=0}^\infty \left(\alpha_{m,n:i} a_{(-m-n-i)}b_{(i)}u+\alpha_{m,n;i}b_{(-m-n-i)}a_{(i)}u\right)
\end{align*}
for simplicity.
Then by \eqref{aiuewy}, we have $(a_{(-m)}b_{(-n)}\1)_{(-1)}u=a_{(-m)}b_{(-n)}u+w$.
Hence
\begin{align*}
\overline{\eta}(a_{(-m)}b_{(-n)}\1)\cdot \eta(u)&=\overline{\eta}((a_{(-m)}b_{(-n)}\1)_{(-1)}u)+\overline{\phi}(a_{(-m)}b_{(-n)}\1,u)\\
&=\overline{\eta}(a_{(-m)}b_{(-n)}u)+\overline{\eta}(w)+\overline{\phi}(a_{(-m)}b_{(-n)}\1,u).
\end{align*}
By Lemma \ref{apjrj}, we see that
\[
\overline{\phi}(a_{(-m)}b_{(-n)}\1,u)=\overline{\phi}(\1,b_{(-n)}a_{(-m)}u)=\overline{\eta}(b_{(-n)}a_{(-m)}u)
\]
since $m,n\geq 2$.
Now Commutativity formula \eqref{comm1} shows the lemma.
\end{proof}
Lemma \ref{uaygsdf} implies that if $m,n\geq 2$ then  $\overline{\eta}(a_{(-m)}b_{(-n)}\1)\cdot \overline{\eta}(u)$ is in $\image\overline{\eta}$ for any $a,b,u\in V$.
We also find that by Proposition \ref{cdyuga},  
\begin{align}\label{wjhsbd}
\overline{\eta}(a_{(-m)}b_{(-n)}u)\equiv \frac{1}{2}\overline{\eta}(a_{(-m)}b_{(-n)}\1)\cdot \overline{\eta}(u)
\end{align}
modulo $\overline{\eta}(\mathcal{F}_{\wt{a}+\wt{b}+k-1}(V))$ for any $m,n\geq 2$, homogeneous $a,b\in V$ and $u\in\mathcal{F}_k(V)$.

Let $a^i\in V$, $m_i\in\Zplus$ $(i=1,2,\cdots,r)$.
We introduce a notation
\[
\Phi(a^1,\cdots,a^r)_{m_1,\cdots,m_r}:=(m_1-1)!\cdots(m_r-1)!{a^1}_{(-m_1)}\cdots {a^r}_{(-m_r)}\1.
\]
It is clear that
\[
L_{-1}\Phi(a^1,\cdots,a^r)_{m_1,\cdots,m_r}=\sum_{i=1}^r\Phi(a^1,\cdots,a^r)_{m_1,\cdots,m_i+1,\cdots,m_r}
\]
holds.
Thus \eqref{iwweyyo} shows that
\begin{align}\label{asidfy}
\sum_{i=1}^r\overline{\eta}\left(\Phi(a^1,\cdots,a^r)_{m_1,\cdots,m_i+1,\cdots,m_r}\right)=0.
\end{align}
By using this identity for $r=2$, we have the following lemma.
\begin{lemma}\label{suidygf}
For $a,b\in V$, $m,n\in\Zplus$ and $-m<l< n$,
\begin{align}\label{uasydfg}
\overline{\eta}\left(\Phi(a,b)_{m,n}\right)=(-1)^{l}\overline{\eta}\left(\Phi(a,b)_{m+l,n-l}\right).
\end{align}
In particular, if $m+n$ is odd, then $\overline{\eta}\left(\Phi(a,a)_{m,n}\right)=0$.
\end{lemma}
\begin{proof}
By \eqref{asidfy}, we have $\overline{\eta}\left(\Phi(a,b)_{m,n}\right)=-\overline{\eta}\left(\Phi(a,b)_{m+1,n-1}\right)$ for $m\geq 1$ and $n\geq 2$.
Induction on $l>0$ proves \eqref{uasydfg}.
We can also show that \eqref{uasydfg} is valid if $l\leq 0$.

If $m+n$ is odd and $m+n=2k+1$ $(k\in\Zplus)$, then
\[
\overline{\eta}\left(\Phi(a,a)_{m,n}\right)=(-1)^{n-k-1}\overline{\eta}\left(\Phi(a,a)_{k,k+1}\right)=(-1)^{n-k}\overline{\eta}\left(\Phi(a,a)_{k+1,k}\right).
\]
Since $\Phi(a,a)_{k,k+1}\equiv \Phi(a,a)_{k+1,k}$ modulo $L_{-1}V$ by Commutativity formula \eqref{comm1}, we have $\overline{\eta}(\Phi(a,a)_{k+1,k})=0$ and $\overline{\eta}\left(\Phi(a,a)_{m,n}\right)=0$.
\end{proof}


We shall start a proof of Lemma \ref{ashgf}.
For simplicity, we set
\[
{U}_k:=\overline{\eta}(\mathcal{F}_{k}(V))
\]
for any $k\in\Zpos$.
For $s\geq1$, we set
\begin{align*}
\mathcal{B}_{k}^{s;r}&=\haru{\overline{\eta}\left(\Phi(x^1,\cdots,x^s)_{m_1,\cdots,m_r,1^s}\right)}{x^i\in T,\sum\wt{x^i}=k,m_i\geq 2},
\end{align*}
for $k\in\Zpos$ and $r,s\in\Zpos$ with $r\leq s$, where $\Phi(x^1,\cdots,x^s)_{1^{s}}$ denotes $\Phi(x^1,\cdots,x^s)_{1,\cdots,1}$. 
We define $\mathcal{B}_{k}^{s;r}=0$ if $\sum\wt{x^i}\neq k$ for any $x^i\in T$.
We note that if $s>k$ then $\mathcal{B}_{k}^{s;r}=0$ because $T$ contains only homogeneous vectors of positive weight.
We also find that each $\mathcal{B}_{k}^{s,0}$ is finite dimensional for each $k,s$.
\begin{lemma}\label{claim1}
 $(\mathcal{B}_{k}^{s;r}+{U}_{k-1})/{U}_{k-1}$ is finite dimensional for any $r,k\geq 1$.
\end{lemma}
\begin{proof}
In the case $s=1$, $\mathcal{B}_{k}^{1;1}=\haru{\overline{\eta}(\Phi(x)_{m})}{x\in  T,\wt{x}=k,m\geq 2}$.
Since $\Phi(x)_{m}\in L_{-1}V$ for $m\geq 2$, we have $\mathcal{B}_{k}^{1;1}=0$.

Next we consider the case $s\geq 2$.
If $r\geq 2$, then \eqref{wjhsbd} shows that
\begin{align*}
&\overline{\eta}(\Phi(x^1,\cdots,x^s)_{m_1,\cdots,m_r,1^{s-r}})\\
&\equiv \frac{1}{2}\overline{\eta}(\Phi(x^1,x^2)_{m_1,m_2})\cdot\overline{\eta}(\Phi(x^3,\cdots,x^s)_{m_3.\cdots,m_r,1^{s-r}})
\end{align*}
modulo $U_{k-1}$.
It follows from Lemma \ref{suidygf} that
\[
\overline{\eta}(\Phi(x^1,x^2)_{m_1,m_2})=(-1)^{m_2}\overline{\eta}(\Phi(x^1,x^2)_{m_1+m_2-2,2}).
\]
Therefore we have   
\begin{align*}
\overline{\eta}(\Phi(x^1,\cdots,x^s)_{m_1,\cdots,m_r,1^{s-r}})\equiv (-1)^{m_2}\overline{\eta}(\Phi(x^1,\cdots,x^s)_{m_1+m_2-2,2,m_3,\cdots,m_r,1^{s-r}})
\end{align*}
modulo $U_{k-1}$.
Proposition \ref{cdyuga} (1) and the argument above show that
\begin{align*}
&\overline{\eta}(\Phi(x^1,\cdots,x^s)_{m_1,\cdots,m_r,1^{s-r}})\\
&\equiv (-1)^{\sum_{i=2}^{r}m_i}\overline{\eta}(\Phi(x^1,\cdots,x^s)_{2+\sum_{i=1}^{r}(m_i-2),2,\cdots,2,1^{s-r}})
\end{align*}
modulo $U_{k-1}$.
Hence we have
\begin{align*}
&\mathcal{B}_{k}^{s;r}+U_{k-1}\\
&=\haru{\overline{\eta}(\Phi(x^1,\cdots,x^s)_{m,2,\cdots,2,1^{s-r}})}{x^i\in  T,\sum\wt{x^i}=k,m\geq 2}+U_{k-1}\\
&=\sum_{l=2}^k \mathcal{B}_{l}^{2;2}\cdot \haru{\overline{\eta}(\Phi(x^3,\cdots,x^{s})_{2,\cdots,2,1^{s-r}})}{x^i\in T,\sum\wt{x^i}=k-l}+U_{k-1}.
\end{align*}
From the assumption of Lemma \ref{ashgf}, $\mathcal{B}_{l}^{2;2}$ is finite dimensional for any $l$.
We also find that the subspace spanned by $\overline{\eta}(\Phi(x^3,\cdots,x^{s})_{2,\cdots,2,1^{r-s}})$ with $x^i\in T$ satisfying $\sum\wt{x^i}=k-l$ is also finite dimensional for each $l$.
Hence $(\mathcal{B}_{k}^{s;r}+U_{k-1})/U_{k-1}$ is finite dimensional if $r\geq 2$.

Finally we consider the case $r=1$. 
For $m\geq 3$, by \eqref{asidfy} and Proposition \ref{cdyuga} (1), we have 
\begin{align*}
\overline{\eta}(\Phi(x^1,\cdots,x^{s})_{m,1^{s-1}})=-\sum_{i=2}^{s}\overline{\eta}(\Phi(x^1,\cdots,x^{s})_{m-1,1^{i-2},{2},1^{s-i}})\in \mathcal{B}_{k}^{s;2}+U_{k-1}
\end{align*}
for $x^i\in T$ with $\sum\wt{x^i}=k$. 
This shows that 
\[
\mathcal{B}_{k}^{s;1}\subset \haru{\overline{\eta}(\Phi(x^1,\cdots,x^{s})_{2,1,\cdots,1})}{x^i\in T}+ \mathcal{B}_{k}^{s;2}+U_{k-1}.
\]
Hence $(\mathcal{B}_{k}^{s;1}+U_{k-1})/U_{k-1}$ is finite dimensional. 
This completes the proof of Lemma \ref{claim1}.
\end{proof}
\begin{lemma}\label{claim2}
${U}_k$ is finite dimensional for any $k\geq 0$.
\end{lemma}
\begin{proof}
By definition and Proposition \ref{cdyuga} (3), for $k\geq 1$,   
\begin{align}\label{iajdsfh}
U_{k}=\sum_{s=1}^{k}\sum_{r=0}^{s} \mathcal{B}_{k}^{s;r}+U_{k-1}. 
\end{align}
Since $\mathcal{B}_{k}^{s,0}$ is finite dimensional, Lemma \ref{claim1} proves the quotient space $U_k/U_{k-1}$ is of finite dimension.
Thus $U_0=\C\1$ implies that $U_k$ is finite dimensional for all $k\in\Zpos$.
\end{proof}

Now we suppose that $V$ is $C_2$-cofinite and show that $U_{k}=U_{k-1}$ for sufficiently large $k$.
Since $V$ is $C_2$-cofinite, there exists $p\in\Zplus$ such that $\Phi(x^1,\cdots,x^t)_{1^t}$ is in $\mathcal{F}_{\sum \wt{x^i}-1}(V)$ for any $t\geq p$ and $x^i\in T$.
\begin{lemma}\label{siudfhiw}
For any $k\geq s\geq r\geq 1$, if $s-r\geq p$ then $\mathcal{B}_{k}^{s,r}\subset U_{k-1}$. 
\end{lemma}
\begin{proof}
If $s-r\geq p$, then we have $\Phi(x^1,\cdots,x^s)_{m_1,\cdots,m_r,1^{s-r}}\in\mathcal{F}_{\sum \wt{x^i}-1}(V)$ for $m_i\in\Zplus$.
Thus $\mathcal{B}_{k}^{s,r}\subset U_{k-1}$.  
\end{proof}

We set $\kappa:=\max\{\wt{x}|x\in T\}$ and $\nu:=\sharp{T}$.
One notices that $\mathcal{B}_{k}^{s,r}=0$ if $s\leq [k/\kappa]$, where $[x]$ denotes the maximal integer less than or equal to $q$ for $q\in\Q$.
Therefore by \ref{siudfhiw}, if $k\geq p\kappa$ then 
\begin{align}\label{skei}
U_{k}=\sum_{s=[k/\kappa]+1}^k\sum_{r=s-p+1}^{s}\mathcal{B}_{k}^{s,r}+U_{k-1}.
\end{align}

\begin{lemma}\label{claim3} If $k\geq  \kappa(5\nu+p)$, then $\mathcal{B}_{k}^{s;r}\subset U_{k-1}$ for any $s\geq r \geq 1$.
\end{lemma}
\begin{proof}
Since $k/\kappa\geq 5\nu+p$, by \eqref{skei}, we may assume that $s\geq 5\nu+p+1$ and $r\geq 5\nu+1$. 
For such $k,s,r$ we show that the image of $v:=\Phi(x^1,\cdots,x^s)_{m_1,\cdots,m_r,1^{s-r}}$ by $\overline{\eta}$ is in $U_{k-1}$ for any $x^i\in T$ and $m_i\geq 2$. 

Since $r\geq 5\nu+1$, there are $1\leq  i_1<\cdots<i_6\leq r$ such that $x^{i_1}=\cdots=x^{i_6}$. 
Thus by Proposition \ref{cdyuga} (1), we may assume that $x^1=x^2=x^3=x^4=x^5=x^6$, say $y$.
The proof of Lemma \ref{claim1} implies that we  may also assume that $m_2=\cdots=m_r=2$. 

Suppose that $m_1\geq 3$. 
Then we have
\begin{align*}
\overline{\eta}(v)\equiv &\frac{1}{2}\overline{\eta}(\Phi(y,y)_{m_1,2})\cdot\overline{\eta}(\Phi(x^3,\cdots,x^s)_{2,\cdots,2,1^{s-r}})\\
\equiv&-\frac{1}{2}\overline{\eta}(\Phi(y,y)_{m_1-1,3})\cdot\overline{\eta}(\Phi(x^3,\cdots,x^s)_{2,\cdots,2,1^{s-r}})\\
\equiv &-\overline{\eta}(\Phi(x^1,\cdots,x^s)_{m_1-1,3,2,\cdots,2,1^{s-r}})\\
\equiv &-\overline{\eta}(\Phi(x^2,x^3,x^1,x^4,\cdots,x^s)_{3,2,m_1-1,2,\cdots,2,1^{s-r}})\\
\equiv &-\frac{1}{2}\overline{\eta}(\Phi(y,y)_{3,2})\cdot\overline{\eta}(\Phi(y,y,x^5,\cdots,x^s)_{m_1-1,2,\cdots,2,1^{s-r}})\\
=&0
\end{align*}
modulo $U_{k-1}$. 
Hence $\overline{\eta}(v)\in U_{k-1}$. 

We consider the case $m_1=2$. 
In this case, $v={y_{(-2)}}^6w$ with $w=\Phi(x^7,\cdots,x^s)_{2,\cdots,2,1^{s-r}}$.
Then Associativity formula shows that 
\begin{align*}
v&\equiv (y_{(-1)}y)_{(-3)}{y_{(-2)}}^4w-2y_{(-3)}y_{(-1)}{y_{(-2)}}^4w\\
&\equiv  {(y_{(-1)}y)_{(-3)}}^2{y_{(-2)}}^2w-2 (y_{(-1)}y)_{(-3)}y_{(-3)}y_{(-1)}{y_{(-2)}}^2w-2y_{(-3)}y_{(-1)}{y_{(-2)}}^4w\\
&\equiv  {(y_{(-1)}y)_{(-3)}}^3w-2{(y_{(-1)}y)_{(-3)}}^2y_{(-3)}y_{(-1)}w\\
&\quad-2 (y_{(-1)}y)_{(-3)}y_{(-3)}y_{(-1)}{y_{(-2)}}^2w-2y_{(-3)}y_{(-1)}{y_{(-2)}}^4w
\end{align*}
modulo $\mathcal{F}_{k-1}(V)$. 
The images of the last two terms by $\overline{\eta}$  are in $U_{k-1}$ because both of them are congruent to vectors of the form $y_{(-3)}y_{(-2)}u$ for some suitable vectors $u$ modulo $\mathcal{F}_{k-1}(V)$ and $\overline{\eta}(y_{(-3)}y_{(-2)}u)\equiv \frac{1}{2}\overline{\eta}(y_{(-3)}y_{(-2)}\1)\cdot\overline{\eta}(u)=0$ modulo $U_{k-1}$. 

As for the first and second terms, we see that they are congruent to vectors of the form $\Phi(a,a,b,u)_{3,3,3,1}$ for some suitable vectors $a,b,u$ modulo $\mathcal{F}_{k-1}(V)$. 
Then by using the same argument in the proof of Lemma \ref{claim1}, we have 
\begin{align*}
\overline{\eta}(\Phi(a,a,b,u)_{3,3,3,1})&\equiv -\overline{\eta}(\Phi(a,a,b,u)_{3,4,2,1})\\
&\equiv -\frac{1}{2}\overline{\eta}(\Phi(a,a)_{3,4})\cdot\overline{\eta}(\Phi(b,u)_{2,1})=0.
\end{align*}
Therefore we have $\overline{\eta}(v)\in U_{k-1}$. 
Consequently one has $B_{k}^{s,r}\subset U_{k-1}$ for $k\geq  \kappa(5\nu+p)$. 
\end{proof}

Finally we have the following lemma.
\begin{lemma}\label{claim5}
$U_k=U_{k-1}$ for $k\geq  \kappa(5\nu+p)$.
\end{lemma}
\begin{proof}
By \eqref{skei} and Lemma \ref{claim3}, we have $U_{k}=U_{k-1}$ for $k\geq  \kappa(5\nu+p)$.
\end{proof}


Now we can prove Lemma \ref{ashgf} and Theorem \ref{fbeygfw}.

\noindent
{\it  Proof of Lemma \ref{ashgf}.}  
By Lemma \ref{claim5}, $\image \overline{\eta}=U_k$ for sufficiently large $k$.  
Hence by Lemma \ref{claim2} we have Lemma \ref{ashgf}. \hfill{$\Box$}
\vskip1ex

\noindent
{\it Proof of Theorem \ref{fbeygfw}.}
It follows from Corollary \ref{isuhds} and Lemma \ref{ashgf}.
\hfill{$\Box$}
 
\section{$C_2$-cofiniteness of $\widetilde{L(c_{p,q},0)}$}\label{sect5}
In this section we find some identities in $\widetilde{V}/C_2(\widetilde{V})$ concerning with the Virasoro vector and show that $\widetilde{L(c_{p,q},0)}$ is $C_2$-cofinite for any coprime $p,q\geq 2$.  
In this section we call vectors of the form $L_{m_1}\cdots L_{m_r}\1$ with $m_i<0$ a monomial type vector of length $r$.

We set $c_{m,n;i}:=\alpha_{m,n:i}+\alpha_{n,m:i}$ for $m,n\geq 1$.
Then by Lemma \ref{aiuewy},
\begin{align}
\begin{split}\label{ahdfss2}
(L_{-m}L_{-n}\1)_{(-1)}L_{-p}L_{-q}\1&=L_{-m}L_{-n}L_{-p}L_{-q}\1\\
&\quad +\sum_{i=0}^\infty c_{m-1,n-1;i} L_{-m-n+1-i}L_{i-1}L_{-p}L_{-q}\1
\end{split}
\end{align}
holds for $m,n,p,q\geq 2$.
By calculating $L_{i-1}L_{-p}L_{-q}\1$ for $i\geq 0$ and using the fact that $L_{i-1}\1=0$ for $i\geq 0$, \eqref{ahdfss2} implies    
\begin{align}
\begin{split}\label{urjubas}
&(L_{-m}L_{-n}\1)_{(-1)}L_{-p}L_{-q}\1\\
&=L_{-m}L_{-n}L_{-p}L_{-q}\1\\
&\quad+\sum_{i=0}^{p-1} c_{m-1,n-1:i} (p+i-1) L_{-m-n+1-i}L_{-q}L_{i-p-1}\1\\
&\quad+\sum_{i=0}^{q-1} c_{m-1,n-1:i} (q+i-1) L_{-m-n+1-i} L_{-p}L_{i-q-1}\1\\
&\quad+\sum_{i=0}^{p+q-1} c_{m-1,n-1:i} (p+i-1)(q+i-p-1) L_{-m-n+1-i}L_{i-p-q-1}\1\\
&\quad+c_{m-1,n-1;p+1}\binom{q+1}{3}\frac{c_V}{2}L_{-m-n-p}L_{-q}\1\\
&\quad+c_{m-1,n-1;q+1}\binom{p+1}{3}\frac{c_V}{2}L_{-m-n-q}L_{-p}\1\\
&\quad+c_{m-1,n-1;p+q+1} (2p+q)\binom{q+1}{3}\frac{c_V}{2}L_{-m-n-p-q}\1.
\end{split}
\end{align}
By applying $\overline{\eta}$ to the both hand sides in \eqref{urjubas} we shall find some identities between $\overline{\eta}(L_{-m+2}L_{-2}\w)$ and $\overline{\eta}(L_{-m}\w)$ in $\widetilde{V}/(C_2(\widetilde{V})$ for sufficiently large integer $m$.
To get them, we need more lemmas.

First it is clear from \eqref{iwweyyo} that
\begin{align}\label{aisd}
\overline{\eta}(L_{-m}\1)=0 \quad\text{for }m\geq 3.
\end{align}
 
We next have the following lemma by applying Lemma \ref{suidygf} to $a=b=\w$.
\begin{lemma}\label{kjasd}
For $m,n\geq 2$,
\[
\overline{\eta}(L_{-m}L_{-n}\1)=(-1)^{n}\binom{m+n-4}{n-2}\overline{\eta}(L_{-m-n+2}\w).
\]
Moreover if $m+n$ is odd then $\overline{\eta}(L_{-m}L_{-n})=0$.
\end{lemma}

As for monomial type vectors of length $3$, we have the following two lemmas.
\begin{lemma}\label{siuygwe}
If $m,n,l\geq 3$, then
\begin{align*}
\overline{\eta}(L_{-m}L_{-n}L_{-l}\1)=-f(m,n,l)\overline{\eta}(L_{-m-n-l+2}\w),
\end{align*}
where
\begin{align}\label{weoiur}
\begin{split}
f(m,n,l)&=\frac{1}{2}\left((m-n)\binom{m+n+l-4}{l-2}(-1)^l \right.\\
&\qquad+(m-l)\binom{m+n+l-4}{n-2}(-1)^{m+l}\\
&\qquad+\left.(n-l)\binom{m+n+l-4}{m-2}(-1)^m\right).
\end{split}
\end{align}
\end{lemma}
\begin{proof}
It follows from the commutation relations of the Virasoro algebra that  
\begin{align*}
L_{-l}L_{-n}L_{-m}\1&=L_{-m}L_{-n}L_{-l}\1+(n-l)L_{-n-l}L_{-m}\1+(m-l)L_{-n}L_{-l-m}\1\\
&\quad +(m-n)L_{-m-n}L_{-l}\1.
\end{align*}
Since $m,n,l\geq 3$, we have
\begin{align*}
\overline{\eta}(L_{-l}L_{-n}L_{-m}\1)&=\overline{\phi}(L_{-l}L_{-n}L_{-m}\1,\1)\\
&=-\overline{\phi}(\1,L_{-m}L_{-n}L_{-l}\1)\\
&=-\overline{\eta}(L_{-m}L_{-n}L_{-l}\1)
\end{align*}
by Lemma \ref{apjrj}.
It follows from Lemma \ref{kjasd} that
\begin{align}
\overline{\eta}(L_{-m-n}L_{-l})&=\binom{m+n+l-4}{l-2}(-1)^{l}\overline{\eta}(L_{-m-n-l+2}
w),\\
\overline{\eta}(L_{-n}L_{-l-m})&=\binom{m+l+n-4}{n-2}(-1)^{m+l}\overline{\eta}(L_{-m-n-l+2}\w),\\
\overline{\eta}(L_{-l-n}L_{-m})&=\binom{m+n+l-4}{m-2}(-1)^{m}\overline{\eta}(L_{-m-n-l+2}\w).
\end{align}
Therefore,
\begin{align*}
-\overline{\eta}(L_{-m}L_{-n}L_{-l})=\overline{\eta}(L_{-m}L_{-n}L_{-l})+2f(m,n,l)\overline{\eta}(L_{-m-n-l+2}\w)
\end{align*}
by means of $f(m,n,l)$ given in \eqref{weoiur}.
This proves the lemma.
\end{proof}

Next we consider the vector $\overline{\eta}(L_{-m}L_{-n}\w)$ for $m,n\geq 3$.
By Lemma \ref{aiuewy}, for $m,n\geq 2$, one gets the following identity:
\begin{align}\label{qwiueyg}
\begin{split}
(L_{-m}L_{-n}\1)_{(-1)}\w&=L_{-m}L_{-n}\w+c_{m-1,n-1:0} L_{-m-n+1}L_{-3}\1\\
&\quad +2 c_{m-1,n-1:1} L_{-m-n}\w +c_{m-1,n-1;3}\frac{c_V}{2}L_{-m-n-2}\1.
\end{split}
\end{align}
Hence by \eqref{aisd}, \eqref{qwiueyg} and Lemma \ref{kjasd} we have
\begin{align}
\begin{split}\label{iaoiq}
\overline{\eta}((L_{-m}L_{-n}\1)_{(-1)}\w)&=\overline{\eta}(L_{-m}L_{-n}\w)-d_{m-1,n-1}\overline{\eta}(L_{-m-n}\w)
\end{split}
\end{align}
for $m,n\geq 2$, where
\[
d_{m,n}=(m+n)c_{m,n:0}-2 c_{m,n:1}.  
\]
We also note that
\begin{align}\label{asjdf}
\overline{\eta}(L_{-m}L_{-n}\1)\cdot\overline{\eta}(\w)=\overline{\eta}((L_{-m}L_{-n}\1)_{(-1)}\w)+\overline{\phi}(L_{-m}L_{-n}\1,\w)
\end{align}
for $m,n\geq 2$.
Hence if $m,n\geq 3$, then \eqref{iaoiq} and Lemma \ref{apjrj} give
\begin{align*}
\overline{\eta}(L_{-m}L_{-n}\1)\cdot\overline{\eta}(\w)
&=\overline{\eta}(L_{-m}L_{-n}\w)+\overline{\eta}(L_{-n}L_{-m}\w)-d_{m-1,n-1}\overline{\eta}(L_{-m-n}\w)\\
&=2\overline{\eta}(L_{-m}L_{-n}\w)+(m-n-d_{m-1,n-1})\overline{\eta}(L_{-m-n}\w).
\end{align*}
Therefore we have
\begin{align}\label{qiweh}
\begin{split}
2\overline{\eta}(L_{-m}L_{-n}\w)=&\overline{\eta}(L_{-m}L_{-n}\1)\cdot\overline{\eta}(\w)\\
&+(-m+n+d_{m-1,n-1})\overline{\eta}(L_{-m-n}\w)
\end{split}
\end{align}
for $m,n\geq 3$.
Finally by Lemma \ref{kjasd} and \eqref{qiweh}, we see that for $m,n\geq 3$,
\begin{align}
\begin{split}\label{asdiuhfg}
2\overline{\eta}(L_{-m}L_{-n}\w)&=\binom{m+n-4}{n-2}(-1)^n\overline{\eta}(L_{-m-n+2}\w)\cdot\overline{\eta}(\w)\\
&+(-m+n+d_{m-1,n-1}) \overline{\eta}(L_{-m-n}\w).
\end{split}
\end{align}

On the other hand, it also follows from \eqref{asjdf} that for $m\geq 3$,
\[
\overline{\eta}(L_{-m}\w)\cdot\overline{\eta}(\w)=\overline{\eta}((L_{-m}\w)_{(-1)}\w)+\overline{\phi}(L_{-m}\w,\w).
\]
If $m\geq 3$, then $\overline{\phi}(L_{-m}\w,\w)=-\overline{\phi}(\w,L_{-m}\w)=-\overline{\phi}(L_{-m}\w,\w)$, and this shows $\overline{\phi}(L_{-m}\w,\w)=0$.
Thus by \eqref{iaoiq}, we have
\begin{align}\label{pasrf}
\overline{\eta}(L_{-m}\w)\cdot\overline{\eta}(\w)=\overline{\eta}(L_{-m}L_{-2}\w)-d_{m-1,1} \overline{\eta}(L_{-m-2}\w).
\end{align}
Therefore by \eqref{asdiuhfg} and \eqref{pasrf}, we have the following lemma.
\begin{lemma}\label{aeuwygfq}
For $m,n\geq 3$,
\begin{align*}
2\overline{\eta}(L_{-m}L_{-n}\w)=\binom{m+n-4}{n-2}(-1)^n\overline{\eta}(L_{-m-n+2}L_{-2}\w)+g(m,n) \overline{\eta}(L_{-m-n}\w),
\end{align*}
where
\begin{align*}
g(m,n)=-\binom{m+n-4}{n-2}(-1)^n d_{m+n-3,1}-m+n+d_{m-1,n-1}.
\end{align*}
\end{lemma}

Now we return to \eqref{urjubas}.
For $m,n\geq 3$ and $p,q\geq 2$, we have
\begin{align*}
&\overline{\eta}(L_{-m}L_{-n}\1)\cdot\overline{\eta}(L_{-p}L_{-q}\1)\\
&=\overline{\eta}((L_{-m}L_{-n}\1)_{(-1)}L_{-p}L_{-q}\1)+\overline{\phi}(L_{-m}L_{-n}\1,L_{-p}L_{-q}\1)\\
&=\overline{\eta}((L_{-m}L_{-n}\1)_{(-1)}L_{-p}L_{-q}\1)+\overline{\eta}(L_{-n}L_{-m}L_{-p}L_{-q}\1)\\
&=\overline{\eta}((L_{-m}L_{-n}\1)_{(-1)}L_{-p}L_{-q}\1)\\
&\quad +\overline{\eta}(L_{-m}L_{-n}L_{-p}L_{-q}\1)+(m-n)\overline{\eta}(L_{-m-n}L_{-p}L_{-q}\1).
\end{align*}
By \eqref{urjubas}, the right hand side is a sum of $2\overline{\eta}(L_{-m}L_{-n}L_{-p}L_{-q}\1)$ and $\overline{\eta}(w)$ such that $w$ is a linear combination of monomial type vectors whose lengths are $2$ or $3$ and wights are $s:=m+n+p+q$.
Lemmas \ref{kjasd}--\ref{aeuwygfq} show that $\overline{\eta}(w)$  is a linear combination of $\overline{\eta}(L_{-s+2}\w)$ and $\overline{\eta}(L_{-s+4}L_{-2}\w)$.
Therefore we get an identity of the form
\begin{align}
\begin{split}\label{aukerqu}
2\overline{\eta}(L_{-m}L_{-n}L_{-p}L_{-q}\1)=&\overline{\eta}(L_{-m}L_{-n}\1)\cdot\overline{\eta}(L_{-p}L_{-q}\1)\\
&+\alpha \overline{\eta}(L_{-s+4}L_{-2}\w)+\beta \overline{\eta}(L_{-s+2}\w)
\end{split}
\end{align}
for some scalars $\alpha,\beta$ which are able to be calculated explicitly by using Lemmas \ref{kjasd}--\ref{aeuwygfq}.

Now we take $m,n,p,q$ so that $m\geq 14$ and even, $n=13-2k$, $p=3+2k$ for $0\leq k\leq 2$ and $q=2$.
Then we see that $s=m+18$ is an even integer greater than or equal to $32$, and $\overline{\eta}(L_{-m}L_{-n}\1)\cdot\overline{\eta}(L_{-p}L_{-q}\1)$ is zero by Lemma \ref{kjasd}.
Thus we have constants $\alpha_{k},\beta_{k}$ such that  
\begin{align}
\begin{split}\label{aukerqu4}
2\overline{\eta}(L_{-s+18}L_{-13+2k}L_{-3-2k}\w)=\alpha_{k} \overline{\eta}(L_{-s+4}L_{-2}\w)+\beta_{k} \overline{\eta}(L_{-s+2}\w).
\end{split}
\end{align}

On the other hand by calculating \eqref{aukerqu} again after changing $(m,n,p,q)$ to $(m,p,n,q)$, we find constants $\gamma_k,\delta_k\in\C$ such that  
\begin{align}
\begin{split}\label{aukerqu2}
2\overline{\eta}(L_{-s+18}L_{-3-2k}L_{-13+2k}\w)=\gamma_{k} \overline{\eta}(L_{-s+4}L_{-2}\w)+\delta_{k} \overline{\eta}(L_{-s+2}\w),
\end{split}
\end{align}
where we note that $\overline{\eta}(L_{-13+2k}\w)=0$.
We also have
\begin{align*}
L_{-m}L_{-p}L_{-n}L_{-q}\1=L_{-m}L_{-n}L_{-p}L_{-q}\1+(n-p)L_{-m}L_{-n-p}L_{-q}\1.
\end{align*}
Hence by \eqref{aukerqu2} and by Lemmas \ref{siuygwe} and \ref{aeuwygfq}, we have an another identity
\begin{align}
\begin{split}\label{aukerqu3}
2\overline{\eta}(L_{-s+18}L_{-13+2k}L_{-3-2k}\w)=\gamma'_{k} \overline{\eta}(L_{-s+4}L_{-2}\w)+\delta'_{k} \overline{\eta}(L_{-s+2}\w).
\end{split}
\end{align}
Finally we have an identity
\begin{align}\label{ghyeru}
\xi_{k} \overline{\eta}(L_{-s+4}L_{-2}\w)+\zeta_{k} \overline{\eta}(L_{-s+2}\w)=0
\end{align}
where $\xi_{k}=\xi_k(s,c_V)=\alpha_{k}-\gamma'_{k}$ and $\zeta_{k}=\zeta_k(s,c_V)=\beta_{k}-\delta'_{k}$.

In Appendix, we give the explicit forms of the coefficients $\xi_{k}(s,c)$ and $\zeta_{k}(s,c)$ for $k=0,1,2$, $s\geq 32$ and $c\in\C$.
By means of the explicit forms of $\xi_{k}$ and $\zeta_{k}$ for $k=0,1,2$, we can show the following lemma.
\begin{lemma}\label{ikasjdasi}
For any even $s\geq 32$ and $c\in \C$, one of the determinants of matrices $\begin{pmatrix}\xi_{0}(s,c)&\zeta_{0}(s,c)\\ \xi_{1}(s,c)&\zeta_{1}(s,c)\end{pmatrix}$, $\begin{pmatrix}\xi_{1}(s,c)&\zeta_{1}(s,c)\\ \xi_{2}(s,c)&\zeta_{2}(s,c)\end{pmatrix}$ and $\begin{pmatrix}\xi_{2}(s,c)&\zeta_{2}(s,c)\\ \xi_{0}(s,c)&\zeta_{0}(s,c)\end{pmatrix}$ is nonzero.
\end{lemma}
We give a proof of Lemma \ref{ikasjdasi} in Appendix.
This lemma implies that $\overline{\eta}(L_{-s+4}L_{-2}\w)=\overline{\eta}(L_{-s+2}\w)=0$ for any even integer $s\geq 32$ and any central charge $c=c_V\in\C$.
We recall $\overline{\eta}(L_{-n}\w)=0$ if $n$ is odd.
Therefore, we have the following theorem.
\begin{theorem}\label{audfqku}
Let $V$ be a vertex operator algebra.
Then $\overline{\eta}(L_{-n}\w)=0$ for $n\geq 30$.
In particular, the subspace $\haru{\overline{\eta}(L_{-n}\w)}{n\geq 2}$ is finite dimensional.  
\end{theorem}

As an application of Theorems \ref{fbeygfw} and \ref{audfqku}, we consider a $2$-cyclic permutation orbifold models of the Virasoro vertex operator algebras.

Let $V=L(c_{p,q},0)$ be the simple Virasoro vertex operator algebra of central charge $c=c_{p,q}$ for coprime integers $p,q\geq 2$.
It is well known that $V$ is strongly generated by $\w$ and $C_2$-cofinite.
Therefore by Theorems \ref{fbeygfw} and \ref{audfqku}, we have the following theorem.
\begin{theorem}
Let $p,q$ be coprime integers greater than or equal to $2$.
Then the $2$-cyclic permutation orbifold model $\widetilde{L(c_{p,q},0)}$ is $C_2$-cofinite.
\end{theorem}

\section{Appendix}\label{sect6}
In this section we give the explicit forms of $\xi_{k}(m,c)$ and $\zeta_{k}(m,c)$ for $k=0,1,2$ and prove Lemma \ref{ikasjdasi}.
We use Mathematica for the computations in this section.

First we give data to compute $\xi_{0}$ and $\zeta_{0}$.
We take $(m,n,p,q)=(m,13,3,2)$ with even integer $m$.
Then we have  
\begin{align*}
\alpha_{0}=&\frac{1}{14!}(m-13)(m+12)(-60354201600 - 25041744000  m\\
&- 11025031680m^2 +  2218757736 m^3 + 4290676052 m^4 + 2061162870 m^5 \\
&+ 561027415 m^6 + 98527338 m^7 + 11580231 m^8 + 907530 m^9+    45565 m^{10} \\
&+ 1326 m^{11} + 17 m^{12}),\\
\beta_{0}=&\frac{1}{13!\cdot 8}  (m-13) (m+12 )(m+13)(m+14)\\
&\times(-5748019200- 2706163200 m - 1031677920m^2 + 416682968m^3 \\
&\quad + 502648380m^4 + 206064690m^5 + 49210811m^6 + 7683234m^7 \\
&+814359m^8 + 58630m^9 + 2769m^{10} + 78m^{11} + m^{12}) \\
&+\frac{91(241m+738)(m-13)}{2(m+2)(m+3)}\binom{m+14}{16}c.
\end{align*}  
As well, we get
\begin{align*}
\gamma_0'=&\frac{1}{14!}(m+12)(784604620800 + 322088054400  m + 146756039040m^2 \\
&\quad -84093309768m^3 -95650195420m^4 - 30695547818m^5 \\
&\quad - 491574005m^6 +2565009941m^7 + 883762815m^8 + 155173941m^9 \\
&\quad + 16374865m^{10}  + 1047527m^{11} + 37505m^{12} + 577 m^{13}),\\
\delta_{0}'=&-\frac{1}{2\cdot 16!}(m+11)(m+12)(m+14)( 264931430400 + 137634854400 m \\
&\quad+ 43159534560 m^2 - 53918986488 m^3 - 38531775476m^4 \\
&\quad- 7272782558m^5 + 1442243915m^6 + 1002076031m^7 \\
&\quad+ 232782927m^8 + 31029951m^9 + 2612465m^{10} +140117m^{11} \\
&\quad + 4489m^{12} + 67m^{13})\\
&+\frac{33(192721m+2502378)(m-3)}{2(m+2)(m+13)}\binom{m+14}{16}c.
\end{align*}
Finally we get
\begin{align*}
\xi_0=&-\frac{2284800}{(m+1)(m+11)(m+13)}\binom{m+14}{17},\\
\zeta_{0}=&\frac{152320( m^3+30 m^2+437 m+2628)}{(m+1)(m+3)(m+11)(m+13)}\binom{m+14}{17}\\
&-\frac{(3168931 m^3+41305158 m^2-26765899 m-366373782)}{(m+2)(m+3)(m+13)}\binom{m+14}{16}c.
\end{align*}

We can compute the following coefficients by the similar way.   
\begin{align*}
\xi_1=&-\frac{17821440}{(m+5)(m+7)(m+13)}\binom{m+14}{17},\\
\zeta_{1}=&\frac{8}{5\cdot 13!}(13m^6+559 m^5+10914 m^4+ 113042m^3+541013 m^2\\
&\quad +743199 m+141660)(m-2)m(m+2)(m+4)(m+6)(m+7)(m+8)\\
&\times (m+10)(m+12)(m+14)\\
&-\frac{15(1595885 m^3+18046170 m^2-32028341 m-411762282)}{(m+2)(m+5)(m+11)}\binom{m+14}{16}c,
\end{align*}
\begin{align*}
\xi_2=&-\frac{18670080}{(m+3)(m+9)(m+13)}\binom{m+14}{17},\\
\zeta_{2}=&\frac{32}{15!}(143m^6+6435 m^5+129602 m^4+ 1390950 m^3+7149347 m^2\\
&\quad +12054375 m+5315868)(m-2)m(m+2)(m+4)(m+5)(m+6)\\
&\times (m+8)(m+10)(m+12)(m+14)\\
&-\frac{156(152331 m^3+1559062 m^2-4454319 m-55321438)}{(m+2)(m+7)(m+9)}\binom{m+14}{16}c.
\end{align*}

\noindent
{\it Proof of Lemma \ref{ikasjdasi}.}
By using the explicit forms of $\xi_k(m,c)$ and $\zeta_k(m,c)$ for $k=0,1,2$, we can describe the determinants $\begin{vmatrix}\xi_{k}&\zeta_{k}\\ \xi_{k+1}&\zeta_{k+1}\end{vmatrix}$ for $k=0,1,2$ as $p_{k}(m)+q_{k}(m) c$ with polynomials $p_k(m)$ and $q_k(m)$ in $m$, where we identify $\xi_{3}$ and $\zeta_{3}$ with $\xi_{0}$ and $\zeta_{0}$ respectively.
The system of equations $p_{k}(m)+q_{k}(m) c=0$ for all $k=0,1,2$ leads three equations $p_{k}(m)q_{k+1}(m)-p_{k+1}(m)q_{k}(m)=0$ for $k=0,1,2$.
The polynomials $p_{k}(m)q_{k+1}(m)-p_{k+1}(m)q_{k}(m)$ for $k=0,1,2$ are given as products of nonzero constants, powers of factors $(m+r)$ with $-2\leq r\leq 18$ and a polynomial $f(m)$ given by  
\begin{align*}
f(m)&= -5823421556567940  - 13295522326219116  m - 7085484924471269m^2 \\
&\quad- 1746250016719384m^3 - 310878749441408m^4 - 41974581663344m^5 \\
&\quad-4071611633914m^6 - 252490022696m^7 - 6600424292m^8 \\
&\quad + 133103900m^9 + 7930183m^{10}.
\end{align*}
Since $f(m+39)$ is a polynomial in $m$ whose coefficients are all positive, for $m\geq 39$, $f(m)\neq 0$.
We also see that $f(m)\neq 0$ for $32\leq m\leq38$.
Therefore if $m\geq 32$, one of $p_{k}(m)q_{k+1}(m)-p_{k+1}(m)q_{k}(m)$ for $k=0,1,2$ is nonzero.
This implies that one of determinants $p_k(m)+q_k(m)c$ is nonzero for $m\geq 32$ and arbitrary $c\in\C$.
\hfill{$\square$}

\end{document}